\documentclass[10pt,reqno]{amsart}

\usepackage[all]{xy}
\usepackage{a4wide}
\usepackage[mathscr]{eucal}

\newtheorem{prop}{Proposition}[section]
\newtheorem{lem}[prop]{Lemma}
\newtheorem{thm}[prop]{Theorem}
\newtheorem{cor}[prop]{Corollary}

\theoremstyle{remark}
\newtheorem{rem}[prop]{Remark}

\theoremstyle{definition}
\newtheorem{defn}[prop]{Definition}

\numberwithin{equation}{section}
\numberwithin{prop}{section}

\newcommand{\GG}{\mathbb{G}}
\newcommand{\RR}{\mathbb{R}}
\newcommand{\CC}{\mathbb{C}}
\newcommand{\NN}{\mathbb{N}}
\newcommand{\II}{\mathbb{I}}
\newcommand{\UU}{\mathbb{U}}

\newcommand{\eps}{\varepsilon}
\newcommand{\cH}{\mathcal{H}}
\newcommand{\cst}{\mathrm{C^*}}
\newcommand{\wst}{\mathrm{W^*}}
\newcommand{\id}{\mathrm{id}}
\newcommand{\comp}{\!\circ\!}
\newcommand{\tens}{\otimes}
\newcommand{\Wtil}{\widetilde{W}}
\newcommand{\tp}{\xymatrix{*+<.7ex>[o][F-]{\scriptstyle\top}}}
\newcommand{\eq}{\approx}
\newcommand{\cc}{\mathrm{c}}
\newcommand{\coinv}{\kappa}
\renewcommand{\Bar}[1]{\overline{#1}}
\renewcommand{\Hat}[1]{\widehat{#1}}
\newcommand{\is}[2]{\left(#1\,\vline\,#2\right)}
\newcommand{\its}[3]{\left(#1\,\vline\,#2\,\vline\,#3\right)}

\newcommand{\WW}{\mathbb{W}}
\newcommand{\Ahu}{\Hat{A}_{\mathrm{u}}}
\newcommand{\Delhu}{\Hat{\Delta}_{\mathrm{u}}}
\newcommand{\ehu}{\Hat{e}^{\mathrm{u}}}
\newcommand{\kaphu}{\Hat{\kappa}^{\mathrm{u}}}

\newcommand{\tauhu}{\Hat{\tau}^{\mathrm{u}}}
\newcommand{\Rhu}{\Hat{R}^{\mathrm{u}}}

\newcommand{\HHA}{\Hat{\!\!\!\:\Hat{A}}}
\DeclareMathOperator{\M}{M}
\DeclareMathOperator{\B}{B}
\DeclareMathOperator{\K}{\mathscr{K}}
\DeclareMathOperator{\cL}{\mathcal{L}}
\DeclareMathOperator{\Mor}{Mor}
\DeclareMathOperator{\Dom}{Dom}
\DeclareMathOperator{\dom}{D}
\DeclareMathOperator{\Rep}{Rep}
\DeclareMathOperator{\Graph}{Graph}
\DeclareMathOperator{\Hom}{Hom}

\begin{document}


\title{From multiplicative unitaries to quantum groups II}

% \date{December 2, 2006}

\author{Piotr M.~So{\l}tan}
\address{Department of Mathematical Methods in Physics\\
Faculty of Physics\\
Warsaw University}
\email{piotr.soltan@fuw.edu.pl}

\author{Stanis\l{}aw L.~Woronowicz}
\address{Department of Mathematical Methods in Physics\\
Faculty of Physics\\
Warsaw University}
\email{stanislaw.woronowicz@fuw.edu.pl}

\thanks{Research partially supported by KBN grants no.~1P03A03626
and 115/E-343/SPB/6.PRUE/DIE50/2005-2008.}

\begin{abstract}
It is shown that all important features of a $\mathrm{C}^*$-algebraic quantum
group $(A,\Delta)$ defined by a modular multiplicative $W$ depend only on the
pair $(A,\Delta)$ rather than the multiplicative unitary operator $W$. The
proof is based on thorough study of representations of quantum groups. As an
application we present a construction and study properties of the universal
dual of a quantum group defined by a modular multiplicative unitary ---
without assuming existence of Haar weights.
\end{abstract}

\maketitle

\section{Introduction}

Building on the pioneering work of Baaj and Skandalis (\cite{bs})
S.L.~Woronowicz introduced in \cite{mu} the class of \emph{manageable
multiplicative unitary operators.} Such multiplicative unitaries were shown to
give rise to very interesting objects. Every such operator $W$ acting on
$\cH\tens\cH$ (where $\cH$ is some separable Hilbert space) gives rise to a
$\cst$-algebra $A\subset\B(\cH)$ with comultiplication $\Delta$ and a lot of
additional structure (\cite[Theorem 1.5]{mu}). This extra structure comes in
the form of the reduced dual $\Hat{A}$, the position of
$W\in\M\bigl(\Hat{A}\tens{A}\bigr)$, the coinverse $\kappa$, unitary coinverse
$R$ and the scaling group $(\tau_t)_{t\in\RR}$. Moreover $A$ comes naturally
with an embedding into $\B(\cH)$, so it inherits the ultraweak topology from the
latter space. All this structure is defined with direct use of $W$.

The importance of manageability was emphasized with appearance the famous paper
\cite{kv} in which Kustermans and Vaes gave a very satisfactory definition of a
locally compact quantum group. They showed that every such object gives rise to
a manageable multiplicative unitary. In a more recent paper \cite{bsv} it is
shown that the conditions of \emph{regularity} and \emph{semi-regularity} are
not satisfied by multiplicative unitaries related to quantum groups.

Meanwhile, in \cite{azb,axb} (and later in \cite{nazb}) new examples of quantum
groups were constructed using the theory developed in \cite{mu}. Only later, in
\cite{VDhaar} it was shown that they fitted into the scheme of locally compact
quantum groups. Moreover the multiplicative unitaries used to define them were
not manageable, but only \emph{modular}. The difference between the latter
notions is superfluous as was later explained in \cite{mmu}. Still the
possibility that two \emph{different} multiplicative unitary operators gave
rise to quantum groups described by \emph{isomorphic} $\cst$-algebras with
comultiplication (preserved by the
isomorphism) remained unexplored until S.L.~Woronowicz noticed in \cite{haar}
that a formula for a right invariant weight on a quantum group defined by a
modular multiplicative unitary could be expressed by one of the operators
involved in the definition of modularity (cf.~Definition \ref{defmmu}). This
formula might give a weight which is infinite on all non zero positive elements,
but if we choose the multiplicative unitary correctly we may find the Haar
measure for our quantum
group.

This development prompted the following question: if we can use different
multiplicative unitaries to give rise to the same pair $(A,\Delta)$, does the
additional structure on $A$ depend on the choice of the multiplicative unitary?
In this paper we give an answer to this question. The rich structure consisting
of the coinverse, unitary coinverse, scaling group, reduced dual, the reduced
bicharacter and the ultraweak topology on $A$ are determined uniquely by the
pair $(A,\Delta)$ in the sense that they do not depend on the choice of the
multiplicative unitary giving rise to $(A,\Delta)$.

Let us briefly describe the contents of the paper.
In the next section we will recall the definition of a modular multiplicative
unitary and state the most important consequences of the definition. We will
define what we mean by a quantum group and  give a precise formulation of
our main result together with its classical interpretation.

Section \ref{reps} is devoted to developing the representation theory of
quantum groups. This is the main tool in the proof of our main result. We will
define and study strongly continuous representations of a quantum group.
Constructions of direct sums, tensor products and contragradient
representations will be presented. The crucial notions of intertwining
operators, equivalence, quasi equivalence and algebras generated by
representations will be discussed. Section \ref{proofm} contains the proof of
our main theorem. The reasoning is based very firmly on the
facts explained in Section \ref{reps}.

As one application of Theorem \ref{main} we will give,
in Section \ref{uniS}, a detailed account of the construction of the universal
dual of a given quantum group. We will reproduce some of the results of
Kustermans (\cite{univLCQG}) in the more general setting of quantum groups
arising from multiplicative unitaries. Again the main tool will be the theory of
representations of quantum groups developed in Section \ref{reps}. The notion
of a universal quantum group $\cst$-algebra will be introduced and properties
of this object will be studied.

Throughout the paper we will freely use the language of $\cst$-algebras
developed for use in the theory of quantum groups. We refer the reader to
papers \cite{lance,unbo,gen} for notions of multiplier algebras, morphisms of
$\cst$-algebras, $\cst$-algebras generated by quantum families of multipliers,
etc.

\section{Definitions and results}\label{intro}

Let us recall the definition of a modular multiplicative unitary. We
shall use the complex conjugate Hilbert space $\Bar{\cH}$ of a given Hilbert
space $\cH$. Its precise definition is given in Subsection \ref{oper}.

\begin{defn}[{\cite[Definition 2.1]{mmu}}]\label{defmmu}
Let $\cH$ be a Hilbert space. A unitary operator $W\in\B(\cH\tens\cH)$ is a
\emph{modular multiplicative unitary} if it is a multiplicative unitary and
there exist two positive selfadjoint operators $\Hat{Q}$ and $Q$ on $\cH$ with
zero kernels and a unitary operator $\Wtil\in\B\bigl(\Bar{\cH}\tens\cH\bigr)$
such that
\[
W\bigl(\Hat{Q}\tens{Q}\bigr)W^*=\Hat{Q}\tens{Q}
\]
and
\[
\its{x\tens{u}}{W}{z\tens{y}}
=\its{\Bar{z}\tens{Q}u}{\Wtil}{\Bar{x}\tens{Q^{-1}}y}
\]
for all $x,z\in\cH$, $u\in\dom(Q)$ and $y\in\dom\bigl(\Hat{Q}\bigr)$.
\end{defn}

\begin{thm}[\cite{mu,mmu}]\label{mmuThm}
Let $\cH$ be a separable Hilbert space and let $W\in\B(\cH\tens\cH)$ be a
modular multiplicative unitary. Let
\begin{align}
A& =\bigl\{(\omega\tens\id)W:\:\omega\in\B(\cH)_*\bigr\}
^{\scriptscriptstyle\|\cdot\|\text{\rm-closure}},\\
\Hat{A}& =\bigl\{(\id\tens\omega)\bigl(W^*\bigr):\:\omega\in\B(\cH)_*\bigr\}
^{\scriptscriptstyle\|\cdot\|\text{\rm-closure}}.
\end{align}
Then
\begin{enumerate}
\item $A$ and $\Hat{A}$ are non degenerate $\cst$-subalgebras of
$\B(\cH)$;\label{mmuThm1}
\item $W\in\M\bigl(\Hat{A}\tens{A}\bigr)$;\label{mmuThm2}
\item\label{mmuThm3} there exists a unique $\Delta\in\Mor(A,A\tens{A})$ such
that
\[
(\id\tens\Delta)W=W_{12}W_{13};
\]
moreover $\Delta$ is coassociative and the sets
\[
\bigl\{(a\tens{I})\Delta(b):\:a,b\in{A}\bigr\}\quad\text{and}\quad
\bigl\{\Delta(a)(I\tens{b}):\:a,b\in{A}\bigr\}
\]
are linearly dense subsets of $A\tens{A}$;
\item\label{mmuThm4} there exist a unique closed linear operator $\coinv$ on
the Banach space $A$ such that the set
$\bigl\{(\omega\tens\id)W:\:\omega\in\B(\cH)_*\bigr\}$ is a core for $\coinv$
and
\[
\coinv\bigl((\omega\tens\id)W\bigr)=(\omega\tens\id)\bigl(W^*\bigr);
\]
furthermore for any $a,b\in\Dom(\coinv)$ the product $ab\in\Dom(\coinv)$ and
$\coinv(ab)=\coinv(b)\coinv(a)$, the image of $\coinv$ coincides with
$\Dom(\coinv)^*$ and $\coinv\bigl(\coinv(a)^*\bigr)^*=a$ for any
$a\in\Dom(\coinv)$;
\item\label{kap} there exists a unique one parameter group
$(\tau_t)_{t\in\RR}$
of $*$-automorphisms of $A$ and a unique  ultraweakly continuous involutive
$*$-anti-automorphism $R$ of $A$ such that $R\comp\tau_t=\tau_t\comp R$ for all
$t\in\RR$ and $\coinv=R\comp\tau_{i/2}$;
\end{enumerate}
\end{thm}

The objects $\coinv$, $(\tau_t)$ and $R$ appearing in Statement \eqref{kap} of
Theorem \ref{mmuThm} are referred to as the \emph{coinverse}, \emph{scaling
group} and \emph{unitary coinverse}.

All results of the fundamental paper \cite{mu} have been formulated for
multiplicative unitaries acting on separable Hilbert spaces. For this reason we
shall restrict our attention solely to such spaces. In other words, from now on
all Hilbert spaces are assumed to be separable. Moreover if existence of a
certain Hilbert space is a part of statement of a theorem (see e.g.~Proposition
\ref{later}) then it can be shown that this Hilbert space is (or can be chosen)
separable. Still, many of our results are also true if the Hilbert spaces are
of arbitrary dimension.

We shall consider pairs $(A,\Delta)$ consisting of a $\cst$-algebra $A$ and a
morphism $\Delta\in\Mor(A,A\tens{A})$. We say that two such pairs $(A,\Delta_A)$
and $(B,\Delta_B)$ are isomorphic if there is an isomorphism $\Phi\in\Mor(A,B)$
such that
\begin{equation}\label{QGMor}
\Delta_B\comp\Phi=(\Phi\tens\Phi)\comp\Delta_A.
\end{equation}

\begin{defn}\label{DefQG}
Let $A$ be a $\cst$-algebra and $\Delta\in\Mor(A,A\tens{A})$. We say that a
the pair $\GG=(A,\Delta)$ is a \emph{quantum group} if there exists a modular
multiplicative unitary such that $(A,\Delta)$ is isomorphic to the
$\cst$-algebra with comultiplication associated to $W$ in the way described
in Theorem \ref{mmuThm}. In such a case we shall say that $W$ is a modular
multiplicative unitary \emph{giving rise} to the quantum group $\GG$.
\end{defn}

The following definition has been proposed e.g.~in \cite[Page 237]{pu-slw}.

\begin{rem}\label{manag}
Let us note that the results of \cite{mmu} guarantee that $\GG=(A,\Delta)$ is a
quantum group if and only if there exists a \emph{manageable} multiplicative
unitary (\cite[Definition 1.2]{mu}) giving rise to $\GG$.
\end{rem}

The aim of this paper is to provide justification for Definition
\ref{DefQG}.

A modular multiplicative unitary $W$ on a Hilbert space $\cH$ gives
rise to a quantum group $\GG=(A,\Delta)$ as described in Theorem \ref{mmuThm},
but it also produces another quantum group called the \emph{reduced dual} of
$\GG$. This is the quantum group $\Hat{\GG}=\bigl(\Hat{A},\Hat{\Delta}\bigr)$,
where $\Hat{A}$ is the $\cst$-subalgebra of $\B(\cH)$ described in Theorem
\ref{mmuThm} and $\Hat{\Delta}$ is given by
\begin{equation}\label{DAh}
\Hat{\Delta}(x)=\sigma\bigl(W^*(I\tens{x})W\bigr),
\end{equation}
where $\sigma$ is the flip on the tensor product $\Hat{A}\tens\Hat{A}$.
One possible modular multiplicative unitary giving rise to $\Hat{\GG}$ is
$\Hat{W}=\Sigma W^*\Sigma$, where $\Sigma$ is the flip on $\cH\tens\cH$. It is
clear that the reduced dual of $\Hat{\GG}$ is isomorphic to $\GG$.
Note that equation \eqref{DAh} and Theorem \ref{mmuThm} \eqref{mmuThm3} show
that the element $W\in\M\bigl(\Hat{A}\tens{A}\bigr)$
is a \emph{bicharacter}, i.e.~it satisfies
\begin{equation}\label{known}
(\id\tens\Delta)W=W_{12}W_{13}\qquad\text{and}\qquad
\bigl(\Hat{\Delta}\tens\id\bigr)W=W_{23}W_{13}.
\end{equation}
The quantum group $\Hat{\GG}$ and the bicharacter
$W\in\M\bigl(\Hat{A}\tens{A}\bigr)$ are \emph{a priori} defined in terms of the
modular multiplicative unitary which gives rise to $\GG$, rather than $\GG$
itself.

The elements of the polar decomposition of the coinverse $\coinv$ are
also determined \emph{a priori} by the multiplicative unitary. For example the
scaling group $(\tau_t)$ is given by
\[
\tau_t(a)=Q^{2it}aQ^{-2it},
\]
where $Q$ is one of the operators involved in the definition of modularity. In
addition any modular multiplicative unitary $W\in\B(\cH\tens\cH)$ giving rise
to $\GG$ provides a topology on $A$, namely the restriction of the ultraweak
topology from $\B(\cH)$ to $A$. The anti-automorphism $R$ and the automorphisms
$(\tau_t)_{t\in\RR}$ (also defined through $W$) are continuous for this
topology.

We shall prove that all necessary data of a quantum group $\GG$ are independent
of the choice of modular multiplicative unitary giving rise to $\GG$.

\begin{thm}\label{main}
Let $\GG=(A,\Delta)$ be a quantum group. Choose a
Hilbert space $\cH$ and
a modular multiplicative unitary $W\in\B(\cH\tens\cH)$ giving rise to $\GG$ and
use Theorem \ref{mmuThm} to construct the embedding $A\subset\B(\cH)$ and the
objects $\Hat{A}$, $\coinv$, $R$ and $\left(\tau_t\right)_{t\in\RR}$. Let
$\Hat{\Delta}$ be the comultiplication on $\Hat{A}$ given by \eqref{DAh}.
Then
\begin{enumerate}
\item\label{main1} The ultraweak topology on $A$ inherited from $\B(\cH)$ is
independent of the choice of $\cH$ and $W$.
\item\label{main2} The coinverse $\coinv$, its domain and all elements of polar
decomposition are independent of the choice of $\cH$ and $W$;
\item\label{main3} The reduced dual
$\Hat{\GG}=\bigl(\Hat{A},\Hat{\Delta}\bigr)$ and the
bicharacter $W\in\M\bigl(\Hat{A}\tens{A}\bigr)$ are defined uniquely $($up to
isomorphism$)$ by $\GG$. They do not depend on the choice of $\cH$ and
$W\in\B(\cH\tens\cH)$.
\end{enumerate}
\end{thm}
\noindent
The proof of Theorem \ref{main} will be given in Section \ref{proofm}.

Note that one way to interpret Theorem \ref{main} is to say that a quantum
group $\GG=(A,\Delta)$ is naturally endowed with an analog of the
\emph{class} of the Haar measure. This is because the ultraweak topology on
$A$ (determined uniquely by $\GG$) fixes the von Neumann algebra $A''$
which is the non commutative analog of $L^\infty(\GG)$. Also the set of all
ultraweakly continuous functionals on $A$ plays the role of $L^1(\GG)$.

\section{Representations of quantum groups}\label{reps}

Throughout this section let us fix a quantum group $\GG=(A,\Delta)$.

\begin{defn}
A \emph{strongly continuous unitary representation} $U$ of $\GG$ acting on a
Hilbert space $H$ is a unitary element $U\in\M\bigl(\K(H)\tens{A}\bigr)$
such that $(\id\tens\Delta)U=U_{12}U_{13}$. The class of all strongly
continuous unitary representations of $\GG$ will be denoted by $\Rep(\GG)$.

We shall use the symbol $H_U$ for the Hilbert space on which the
representation $U$ acts: $U\in\M(\K(H_U)\tens{A})$.
\end{defn}

Recall that $\GG$ is defined as the pair $(A,\Delta)$ arising from some modular
multiplicative unitary $W\in\B(\cH\tens\cH)$ for some Hilbert space
$\cH$ (cf.~Section \ref{intro}). It follows from statements
\eqref{mmuThm1} and \eqref{mmuThm2} of Theorem \ref{mmuThm} that
$W\in\M\bigl(\K(\cH)\tens{A}\bigr)$ and so by the first part of \eqref{known}
$W$ is a strongly continuous unitary representation of $\GG$ on $\cH$.

If $H$ is a Hilbert space then the element
$\II_H=I_H\tens{I_A}\in\M\bigl(\K(H)\tens{A}\bigr)$ clearly is a strongly
continuous unitary representation of $\GG$. Such representations are called
\emph{trivial.}

In what follows we shall most of the time omit the words ``strongly continuous
unitary'' and speak simply about representations of $\GG$.

\subsection{Intertwining operators}

Let $U,V\in\Rep(\GG)$ and let $t\in\B(H_U,H_V)$. We say that $t$
\emph{intertwines} $U$ and $V$ if
\begin{equation}\label{inter}
(t\tens{I})U=V(t\tens{I}).
\end{equation}
The above equation may be understood in several contexts. If the $\cst$-algebra
$A$ is faithfully represented on a Hilbert space $\cH$ then $U$ and $V$ become
elements of $\B(H_U\tens\cH)$ and $\B(H_V\tens\cH)$ respectively. In this
situation \eqref{inter} means that
\[
(t\tens{I_\cH})U=V(t\tens{I_\cH}).
\]
Equivalently $t\in\B(H_U,H_V)$ intertwines $U$ and $V$ if and only if for any
$\omega\in{A^*}$ we have
\[
t(\id\tens\omega)(U)=(\id\tens\omega)(V)t.
\]
Finally we can identify of $\M\bigl(\K(H_U)\tens{A}\bigr)$ and
$\M\bigl(\K(H_U)\tens{A}\bigr)$ with the $\cst$-algebras $\cL(H_U\tens{A})$
and $\cL(U_V\tens{A})$ of adjointable maps of the Hilbert $A$-modules
$H_U\tens{A}$ and $H_V\tens{A}$ respectively (\cite[pages 10 \& 37]{lance}).
Then $t\in\B(H_U,H_V)$ intertwines $U$ and $V$ if and only if
\[
(t\tens{I_A})U=V(t\tens{I_A})
\]
as elements of $\cL(H_U\tens{A},H_V\tens{A})$.

Let $U,V\in\Rep(\GG)$. The set of operators intertwining $U$ and $V$ will be
denoted by $\Hom(U,V)$.

The following properties follow immediately from the definition of $\Hom(U,V)$:
\begin{enumerate}
\item $\Hom(U,V)$ is a subspace of $B(H_U,H_V)$ closed in the weak operator
topology;
\item for any $t\in\Hom(U,V)$ we have $t^*\in\Hom(V,U)$;
\item $I_{H_U}\in\Hom(U,U)$;
\item if $T$ is another representation of $\GG$ then for any $t\in\Hom(U,V)$
and $s\in\Hom(V,T)$ we have $st\in\Hom(U,T)$ and composition of intertwining
operators is bilinear for the vector space structures on $\Hom(U,V)$ and
$\Hom(V,T)$, in particular $\Hom(U,U)$ is a $*$-algebra with unit;
\item for any $t\in\Hom(U,V)$ the composition $t^*t$ is a positive element of
the $*$-algebra $\Hom(U,U)$, i.e.~there exists $x\in\Hom(U,U)$ such that
$t^*t=x^*x$.
\end{enumerate}
All this shows that the class $\Rep(\GG)$ of strongly continuous unitary
representations of $\GG$ with intertwining operators as morphisms forms a
\emph{concrete $\wst$-category} as defined in
\cite[Definitions 1.1 \& 2.1]{Wst}.

\subsection{Equivalence and quasi-equivalence}

\begin{defn}
Let $U$ and $V$ be representations of $\GG$.
\begin{enumerate}
\item $U$ is a \emph{subrepresentation} of $V$ if $\Hom(U,V)$ contains
an isometry;
\item $U$ and $V$ are \emph{equivalent} if $\Hom(U,V)$ contains an invertible
operator;
\item $U$ and $V$ are \emph{disjoint} if $\Hom(U,V)=\{0\}$;
\item $U$ and $V$ are \emph{quasi-equivalent} if no subrepresentation of $U$
is disjoint from $V$ and no subrepresentation of $V$ is disjoint from $U$.
\end{enumerate}
\end{defn}

Let $U,V\in\Rep(\GG)$. We write $U\eq{V}$ if $U$ and
$V$ are equivalent. Clearly ``$\eq$'' is an equivalence relation. One can
show that quasi equivalence is also an equivalence relation
(cf.~e.g.~Proposition \ref{later}).

\begin{rem}
Let $U,V\in\Rep(\GG)$. By \cite[Corollary 2.7]{Wst} $U$ and $V$ are equivalent
if and only if $\Hom(U,V)$ contains a unitary operator. Similarly $U$ and $V$
are equivalent if $\Hom(U,V)$ contains an operator with trivial kernel and
dense range. Moreover equivalence is the same thing as isomorphism in the
$\wst$-category $\Rep(\GG)$.
\end{rem}

\subsection{Operations on representations}\label{oper}

\subsubsection{Direct sums}

Let $(U_\alpha)$ be a family of representations of $\GG$. The $\cst$-algebra
\[
\bigoplus_{\alpha}\bigl(\K(H_{U_\alpha})\tens{A}\bigr)
\]
is contained in
\[
\K\bigl(\mathop{\oplus}\limits_{\alpha}
H_{U_\alpha}\bigr)\tens{A}.
\]
Moreover the inclusion is a morphism of $\cst$-algebras (\cite{pseu}).
Therefore we have
\[
\M\Bigl(\bigoplus_{\alpha}\bigl(\K(H_{U_\alpha})\tens{A}\bigr)\Bigr)
\subset\M\Bigl(
\K\bigl(\mathop{\oplus}\limits_{\alpha}
H_{U_\alpha}\bigr)\tens{A}\Bigr).
\]
The $\cst$-algebra on the left hand side consists of norm bounded families
$(T_\alpha)$ of multipliers of the $\cst$-algebras $\K(H_{U_\alpha})\tens{A}$.
Therefore the family $U=(U_\alpha)$ is a unitary element of
$\M\bigl(\K\bigl(\mathop{\oplus}\limits_{\alpha}
H_{U_\alpha}\bigr)\tens{A}\bigr)$
It is not difficult to see that it is a strongly continuous unitary
representation of $\GG$. So defined $U$ is called the \emph{direct sum} of the
representations $(U_\alpha)$. Of course
$H_U=\mathop{\oplus}\limits_{\alpha}H_{U_\alpha}$. Moreover for any
$\alpha$ the canonical injection $H_\alpha\hookrightarrow{H}$ and
projection $H\to{H_\alpha}$ belong to $\Hom(U_\alpha,U)$ and $\Hom(U,U_\alpha)$
respectively. In particular each $U_\alpha$ is a subrepresentation of $U$.

\begin{rem}\label{remgeq}
Let $(U_\alpha)$ be a family representations of $\GG$ and let $U$ be the
direct sum of $(U_\alpha)$. Then for any $\omega\in{A}_*$ and any $\alpha$
we have
\[
\bigl\|(\id\tens\omega)U\bigr\|\geq\bigl\|(\id\tens\omega)U_\alpha\bigr\|.
\]
\end{rem}

In what follows we shall restrict attention to countable direct sums, so that
our Hilbert spaces remain separable.

\subsubsection{Tensor products}

Let $U,V\in\Rep(\GG)$. Then the element
\[
U\tp{V}=U_{13}V_{23}\in\M\bigl(\K(H_U)\tens\K(H_V)\tens{A}\bigr)=
\M\bigl(\K(H_U\tens{H_V})\tens{A}\bigr)
\]
is a representation of $\GG$. The representation $U\tp{V}$ is the \emph{tensor
product} of representations $U$ and $V$.

If $T$ is another representation of $\GG$ then we have
$(U\tp{V})\tp{T}\eq{U}\tp{(V\tp{T})}$, so the tensor product of any finite
number of representations of $\GG$ is associative up to equivalence.
Note that the operation of taking tensor product is not, in general,
commutative. In the worst case $U\tp{V}$ is not equivalent to $V\tp{U}$.
However, even if $U\tp{V}\eq{V}\tp{U}$, then in general, the flip
$\Sigma:H_U\tens{H_V}\to{H_V}\tens{H_U}$ does not intertwine $U\tp{V}$ with
$V\tp{U}$: $\Sigma\not\in\Hom(U\tp{V},V\tp{U})$.

The operation of taking tensor products endows $\Rep(\GG)$ with the structure
of a monoidal $\wst$-category (\cite[page 39]{TKSUN}).

\subsubsection{Contragradient representations}\label{repc}

Let $H$ be a Hilbert space. The complex conjugate space $\Bar{H}$ is defined as
the set of elements $\Bar{x}$, where $x\in{H}$. The vector space structure on
$\Bar{H}$ is given by $\Bar{x}+\Bar{y}=\Bar{x+y}$ and
$\zeta\Bar{x}=\Bar{\Bar{\zeta}x}$
for $\Bar{x},\Bar{y}\in\Bar{H}$ and $\zeta\in\CC$. The Hilbert space structure
on $\Bar{H}$ is obtained by setting
\[
\is{\Bar{x}}{\Bar{y}}=\is{y}{x},
\]
where on the right hand side we use the scalar product in $H$. We have the
natural operation of \emph{transposition} taking operators on $H$ to operators
on $\Bar{H}$. This operation will be denoted by $m\mapsto{m^\top}$: for any
closed operator $m$ on $K$ the operator $m^\top$ is defined by
\[
\left(\begin{array}
{c}\begin{pmatrix}\Bar{x}\\{\Bar{y}}\end{pmatrix}\in\Graph{m^\top}
      \end{array}\right)
\Longleftrightarrow
\left(\begin{array}{c}
\begin{pmatrix}x\\{y}\end{pmatrix}\in\Graph{m^*}
      \end{array}\right).
\]
When restricted to bounded operators, the transposition becomes an
anti-isomorphism of $\cst$-algebras $\B(H)\to\B\bigl(\Bar{H}\bigr)$.

In what follows we shall denote by $R$ the unitary coinverse of the quantum
group $\GG$. This is the ``unitary'' part of the polar decomposition
 of the coinverse $\kappa$ (cf.~Section \ref{intro}, \cite{mu,mmu}).

\begin{prop}\label{UcProp}
Let $U\in\Rep(\GG)$. Then the element
$U^\cc=U^{\top\tens{R}}\in\M\bigl(\K\bigl(\Bar{H_U}\bigr)\tens{A}\bigr)$ is a
strongly continuous unitary representation of $\GG$ acting on
$H_{U^\cc}=\Bar{H_U}$.
\end{prop}

\begin{proof}
Denote by $\sigma$ the flip map on $A\tens{A}$. Using
\cite[Theorem 2.3(5)]{mmu}
and remembering that $\top$ and $R$ are anti-isomorphisms we obtain
\[
\begin{split}
(\id\tens\Delta)U^\cc
&=(\id\tens\Delta)(\top\tens{R})U=(\top\tens\Delta\comp R)U\\
&=\bigl(\top\tens[\sigma\comp(R\tens{R})\comp\Delta]\bigr) U\\
&=\bigl(\top\tens[\sigma\comp(R\tens{R})]\bigr)(\id\tens\Delta)U\\
&=(\id\tens\sigma)(\top\tens{R}\tens{R})(\id\tens\Delta)U\\
&=(\id\tens\sigma)(\top\tens{R}\tens{R})(U_{12}U_{13})\\
&=(\id\tens\sigma)\bigl(\bigl[(\top\tens{R})U\bigr]_{13}
\bigl[(\top\tens{R})U\bigr]_{12}\bigr)\\
&=(\id\tens\sigma)(U^\cc_{13}U^\cc_{12})=U^\cc_{12}U^\cc_{13}.
\end{split}
\]
\end{proof}

\begin{defn}
Let $U$ be a strongly continuous unitary representation of $\GG$.
The \emph{contragradient representation} of $U$ is the strongly continuous
unitary representation $U^\cc$ of $\GG$ defined in Proposition \ref{UcProp}.
\end{defn}

\begin{rem}
In contrast to existing definitions of contragradient representations in
literature (e.g.~\cite[Section 3]{pseudogr}) we have $(U^\cc)^\cc=U$
for any strongly continuous unitary representation $U$ of $\GG$.
\end{rem}

The operation of taking contragradient representation is well compatible with
tensor products. In fact if $U$ and $V$ are representations of $\GG$ then
identifying $\Bar{H_U\tens{H_V}}$ with $\Bar{H_V}\tens\Bar{H_U}$ via the
unitary map
\[
\Bar{H_V}\tens\Bar{H_U}\ni\Bar{y}\tens\Bar{x}\longmapsto\Bar{x\tens{y}}
\in\Bar{H_U\tens{H_V}}
\]
we have
\begin{equation}\label{UVcc1}
(U\tp{V})^\cc=V^\cc\tp{U}^\cc.
\end{equation}

\subsection{Quasi equivalence and tensor products}

\begin{prop}\label{later}
Let $U$ and $V$ be representations of $\GG$. Then the following conditions
are equivalent:
\begin{enumerate}
\item\label{pierwszy} $U$ and $V$ are quasi equivalent;
\item\label{drugi} There exists a Hilbert space $Z$ such that $\II_Z\tp{U}$
and $\II_Z\tp{V}$ are equivalent.
\end{enumerate}
\end{prop}

Note that this result can be formulated without the notion of tensor product
of representations. This is because tensor product with a trivial
representation is expressible as a direct sum. We omit the proof of this
proposition as it follows the lines of proofs of analogous results for
representations of $\cst$-algebras (\cite{dix}). Moreover, in this paper we
shall exclusively use condition \eqref{drugi} of Proposition \ref{later} as the
definition of quasi equivalence. Equivalence of \eqref{pierwszy} and
\eqref{drugi} will not be used.

\subsection{Algebras generated by representations}\label{Algs}

In this subsection we shall describe the algebras generated by representations
of $\GG$. We shall use the notion of a $\cst$-algebra generated by a quantum
family of affiliated elements (\cite[Definition 4.1]{gen}).

Let $U\in\Rep(\GG)$. Then there exists a unique $\cst$-algebra $B_U$ acting non
degenerately on $H_U$ such that $U\in\M(B_U\tens{A})$ and $B_U$ is generated by
$U$. Indeed, by Remark \ref{manag} there is a Hilbert space $\cH$ and a
manageable multiplicative unitary $W\in\B(\cH\tens\cH)$ giving rise to $\GG$.
Then (cf.~\cite[Theorems 1.6 \& 1.7]{mu}) one may take
\begin{equation}\label{BUDef}
B_U=\bigl\{(\id\tens\omega)(U^*):\:\omega\in\B(\cH)_*\bigr\}
^{\scriptscriptstyle\|\cdot\|\text{\rm-closure}}.
\end{equation}
Uniqueness of $B_U$ follows from the remark after \cite[Definition 4.1]{gen}.
In particular $B_U$ is independent of the multiplicative unitary $W$ and
Hilbert space $\cH$ entering \eqref{BUDef}. Note that $B_U$ is unique not only
as a $\cst$-algebra, but also as a subset of $\B(H_U)$.

The next proposition states some basic facts about $\cst$-algebras generated by
representations. We omit the simple proof.

\begin{prop}\label{UeqV2}
Let $\GG$ be a quantum group and let $U,V\in\Rep(\GG)$. Then
\begin{enumerate}
\item\label{UeqV20} if $U$ is a subrepresentation of $V$ then there exists a
$\Phi\in\Mor(B_V,B_U)$ such that $(\Phi\tens\id)V=U$. This morphism maps $B_V$
onto $B_U$ and is continuous for the ultraweak topologies inherited by $B_V$
and $B_U$ from $\B(H_V)$ and $\B(H_U)$ respectively;
\item\label{UeqV21} if $U\eq{V}$ then there is a spatial isomorphism
$\Phi\in\Mor(B_U,B_V)$ such that $(\Phi\tens\id)U=V$;
\item\label{UeqV22} if $Z$ is a Hilbert space and $V=\II_Z\tp{U}$ then there is
an isomorphism $\Phi\in\Mor(B_U,B_V)$ such that $(\Phi\tens\id)U=V$. Moreover
$\Phi$ is a homeomorphism for the ultraweak topologies inherited by $B_U$ and
$B_V$ from $\B(H_U)$ and $\B(H_V)$ respectively.
\end{enumerate}
\end{prop}

From Proposition \ref{UeqV2} \eqref{UeqV21},\eqref{UeqV22} and Proposition
\ref{later}
we immediately get

\begin{cor}\label{QePhi2}
Let $\GG$ be a quantum group and let $U$, $V$ be representations
of $\GG$. Assume that $U$ and $V$ are quasi equivalent. Then there is an
isomorphism $\Phi\in\Mor(B_U,B_V)$ such that
\[
(\Phi\tens\id)U=V.
\]
Moreover $\Phi$ is a homeomorphism for the ultraweak topologies inherited by
$B_U$ and $B_V$ from $\B(H_U)$ and $\B(H_V)$ respectively.
\end{cor}

At the end of this section let us mention an important result about matrix
elements of representations. In it we shall use the strict closure of the
operator $\kappa$ defined on the strictly dense subset $\Dom(\kappa)$ of
$\M(A)$ (cf.~\cite[page 133]{mu}).

\begin{prop}\label{kapU}
Let $U$ be a representation of $\GG$. Then for any $\eta\in\B(H_U)_*$ the
element $(\eta\tens\id)U$ belongs to the domain of $\kappa$ and we have
\[
\kappa\bigl((\eta\tens\id)U\bigr)=(\eta\tens\id)(U^*).
\]
\end{prop}

This proposition is a direct consequence of
\cite[Theorem 1.7 \& Theorem 1.6(4)]{mu} and the fact that for any quantum
group $\GG$ there is a \emph{manageable} multiplicative unitary giving rise to
$\GG$ (cf.~\cite{mmu}).

\subsection{Absorbing representations}

\begin{defn}
Let $\GG$ be a quantum group and let $U$ be a representation of $\GG$.
\begin{enumerate}
\item $U$ is called \emph{right absorbing} if for any
representation $V$ of $\GG$ we have
\[
V\tp{U}\eq\II_{H_V}\tp{U}.
\]
\item $U$ is called \emph{left absorbing} if for any representation
$V$ of $\GG$ we have
\[
U\tp{V}\eq{U}\tp\II_{H_V}.
\]
\end{enumerate}
\end{defn}

\begin{rem}\label{Labs}
Let $\GG$ be a quantum group. By \eqref{UVcc1} a representation $U$
of $\GG$ is right absorbing if and only if $U^\cc$ is left absorbing.
\end{rem}

\begin{prop}
Let $\GG=(A,\Delta)$ be a quantum group and let $U\in\Rep(\GG)$. Let $\pi$ be
a representation of $A$ on the Hilbert space $H_U$ which is covariant in the
sense that for any $a\in{A}$
\begin{equation}\label{covar}
U(\pi(a)\tens{I})U^*=(\pi\tens\id)\Delta(a).
\end{equation}
Then $U$ is right absorbing.
\end{prop}

\begin{proof}
Let $V$ be a representation of $V$. Applying $(\id\tens\pi\tens\id)$ to
both sides of the equation $(\id\tens\Delta)V=V_{12}V_{13}$ we obtain
\[
U_{23}\bigl[(\id\tens\pi)V\bigr]_{12}U_{23}^*=
\bigl[(\id\tens\pi)V\bigr]_{12}V_{13}.
\]
Therefore
\begin{equation}\label{piV}
\bigl[(\id\tens\pi)V\bigr]_{12}^*U_{23}\bigl[(\id\tens\pi)V\bigr]_{12}=
V_{13}U_{23}.
\end{equation}
The right hand side of \eqref{piV} is by definition equal to $V\tp{U}$ while
the left hand side is equivalent to $U_{23}=\II_{H_V}\tp{U}$. This means that
$U$ is right absorbing.
\end{proof}

Note that if $\cH$ is a Hilbert space and $W\in\B(\cH\tens\cH)$ is a modular
multiplicative unitary giving rise to $\GG$ then $W$, viewed as an element of
$\M(\K(\cH)\tens{A})$, is a representation of $\GG$ and the embedding of $A$
into $B(\cH)$ given by $W$ is a covariant representation of $A$. In
particular we have

\begin{cor}\label{Wabs}
Let $\cH$ be a Hilbert space and let $W\in\B(\cH\tens\cH)$ be a modular
multiplicative unitary giving rise to $\GG$. Then the representation
$W\in\M(\K(\cH)\tens{A})$ of $\GG$ on $\cH$ is right absorbing.
\end{cor}

\begin{prop}\label{any2}
Let $\GG$ be a quantum group. Then any two right absorbing representations are
quasi equivalent.
\end{prop}

\begin{proof}
Let $U$ and $V$ be right absorbing representations of $\GG$ and let
$T$ be a left absorbing representation of $\GG$ (one can take e.g.~$T=U^\cc$,
cf.~Remark \ref{Labs}). We have
\[
\II_{H_T}\tp{V}\eq
T\tp{V}\eq
T\tp\II_{H_V}
\]
and
\[
T\tp\II_{H_U}\eq
T\tp{U}\eq
\II_{H_T}\tp{U}.
\]
Clearly $T\tp\II_{H_U}$ and $T\tp\II_{H_V}$ are quasi equivalent.
\end{proof}

\section{Proof of main theorem}\label{proofm}

Let $\GG=(A,\Delta)$ be a quantum group. Let $U$ be a right absorbing
representation of $\GG$. Then There is a unique comultiplication
$\Hat{\Delta}_U$ on $B_U$ such that
\begin{equation}\label{U23U13}
\bigl(\Hat{\Delta}_U\tens\id\bigr)U=U_{23}U_{13}.
\end{equation}
To see this let $\cH$ be a Hilbert space and let $W\in\B(\cH\tens\cH)$ be a
modular multiplicative unitary giving rise to $\GG$. the second part of
\eqref{known} tells us that
\[
\bigl(\Hat{\Delta}\tens\id\bigr)W=W_{23}W_{13}.
\]
Now both $U$ and $W\in\M\bigl(\K(\cH)\tens{A}\bigr)$ are right absorbing
representations of $\GG$ (by Corollary \ref{Wabs}) and by Proposition
\ref{any2} and Corollary \ref{QePhi2} there is an isomorphism
$\Phi\in\Mor\bigl(B_U,\Hat{A}\bigr)$ which is a homeomorphism for the
ultraweak topologies on $B_U\subset\B(H_U)$ and $\Hat{A}\subset\B(\cH)$ and
\begin{equation}\label{PhiUW}
(\Phi\tens\id)U=W.
\end{equation}
Therefore setting
$\Hat{\Delta}_U=(\Phi\tens\Phi)^{-1}\comp\Hat{\Delta}\comp\Phi$ we obtain a
comultiplication on $B_U$ satisfying \eqref{U23U13}.

From what we have seen so far it is clear that $\bigl(B_U,\Hat{\Delta}_U\bigr)$
is a quantum group isomorphic to $\Hat{\GG}=\bigl(\Hat{A},\Hat{\Delta}\bigr)$,
i.e.~the reduced dual of $\GG$ defined by $W$.

Now $U$ could have been any other modular multiplicative unitary giving rise to
$\GG$. It follows that the reduced dual $\bigl(\Hat{A},\Hat{\Delta}\bigr)$ is
independent of the multiplicative unitary giving rise to $\GG$. Moreover, the
ultraweak topology on $\Hat{A}$ is independent of $W$.

Repeating the above reasoning for the quantum group $\Hat{\GG}$ we see that
the ultraweak topology on $A$ (which is $\,\HHA$) is independent
of the
modular multiplicative unitary giving rise to $G$. This proves Statement
\eqref{main1} of Theorem \ref{main}.

We have already shown that the reduced dual $\Hat{\GG}$ is independent of the
choice of modular multiplicative unitary giving rise to $\GG$. The position of
$W$ in $\M\bigl(\Hat{A}\tens{A}\bigr)$ is also fixed uniquely. Indeed, for any
right absorbing representation $U$ we have the isomorphism
$\Phi\in\Mor\bigl(B_U,\Hat{A}\bigr)$ satisfying \eqref{PhiUW}. This proves
Statement \eqref{main3} of Theorem \ref{main}.

Statement \eqref{main2} of Theorem \ref{main} follows from Statements
\eqref{main1} and \eqref{main3}. To see this notice that given a modular
multiplicative unitary $W\in\B(\cH\tens\cH)$, the core of $\kappa$ is determined
by the ultraweak topology inherited by $\Hat{A}$ from $\B(\cH)$. Now this
topology is independent of $W$ while the action of $\kappa$ on this core
depends only on the position of $W$ in $\M\bigl(\Hat{A}\tens{A}\bigr)$
(cf.~Theorem \ref{mmuThm} \eqref{mmuThm4}). The uniqueness of the polar
decomposition of $\kappa$ guarantees that $R$ and $(\tau_t)_{t\in\RR}$ are
independent of the choice of modular multiplicative unitary giving rise to
$\GG$.

From now on we shall write $A_*$ for the space of functionals on $A$ continuous
for the ultraweak topology on $A$ coming from representation of $A$ defined by
any modular multiplicative unitary. The image of any right absorbing
representation in $\M\bigl(\Hat{A}\tens{A}\bigr)$ will be called \emph{the
reduced bicharacter} for $\bigl(\GG,\Hat{\GG}\bigr)$. In what follows the
reduced bicharacter will be denoted by the letter $W$. By a \emph{realization}
of $W$ on a Hilbert space $\cH$ we shall mean any modular multiplicative
unitary acting on $\cH\tens\cH$ giving rise to $\GG$.

\section{Universal dual of a quantum group}\label{uniS}

The aim of this section is to define and analyze the universal dual object of a
given quantum group $\GG=(A,\Delta)$. Such objects were already considered in
\cite[Section 3]{pseu} under the name ``Pontryagin dual''.

\subsection{Maximal representations and universal $\cst$-algebra}

\begin{prop}\label{exMax}
There exists a strongly continuous representation $\WW$ of $\GG$ such that for
any $U\in\Rep(\GG)$ and any $\omega\in{A}_*$ we have
\begin{equation}\label{geq}
\bigl\|(\id\tens\omega)\WW\bigr\|\geq\bigl\|(\id\tens\omega)U\bigl\|.
\end{equation}
\end{prop}

\begin{proof}
Let $(\omega_n)_{n\in\NN}$ be a sequence of elements of $A_*$ which is dense in
$A_*$ and each element is repeated infinitely many times. For any $n\in\NN$
there exists a representation $U_n$ of $\GG$ such that
\[
\bigl\|(\id\tens\omega_n)U_n\bigr\|
\geq\sup_U\bigl\|(\id\tens\omega_n)U\bigr\|-\textstyle{\frac{1}{n}}
\]
where the supremum is taken over all strongly continuous unitary
representations of $\GG$.

We define $\WW$ to be the direct sum of $(U_n)_{n\in\NN}$. Formula \eqref{geq}
follows immediately from the definition of $\WW$. Indeed, given a
representation $U$ of $\GG$, $\omega\in{A}_*$ and $\eps>0$ we can find $n$ such
that $\|\omega-\omega_n\|<\frac{\eps}{3}$ and $n>\frac{3}{\eps}$.
Then
\[
\bigl\|(\id\tens\omega_n)U_n\bigr\|\geq
\bigl\|(\id\tens\omega_n)U\bigr\|-\textstyle{\frac{\eps}{3}}
\]
and
\[
\begin{split}
\bigl|\,\bigl\|(\id\tens\omega)\WW\bigr\|
-\bigl\|(\id\tens\omega_n)\WW\bigr\|\,\bigr|
&\leq\bigl\|(\id\tens\omega)\WW-(\id\tens\omega_n)\WW\bigr\|
<\textstyle{\frac{\eps}{3}},\\
\bigl|\,\bigl\|(\id\tens\omega)U\bigr\|
-\bigl\|(\id\tens\omega_n)U\bigr\|\,\bigr|
&\leq\bigl\|(\id\tens\omega_n)U-(\id\tens\omega)U\bigr\|
<\textstyle{\frac{\eps}{3}}.\\
\end{split}
\]
Therefore
\[
\begin{split}
\bigl\|(\id\tens\omega)\WW\bigr\|&\geq\bigl\|(\id\tens\omega_n)\WW\bigr\|
-\textstyle{\frac{\eps}{3}}\geq\bigl\|(\id\tens\omega_n)U_n\bigr\|
-\textstyle{\frac{\eps}{3}}\\
&\geq\bigl\|(\id\tens\omega_n)U\bigr\|-\textstyle{\frac{2\eps}{3}}
\geq\bigl\|(\id\tens\omega)U\bigr\|-\eps
\end{split}
\]
(cf.~Remark \ref{remgeq}).
\end{proof}

\begin{defn}
Let $\GG$ be a quantum group. A representation $\WW$ fulfilling the
condition of Proposition \ref{exMax} is called \emph{maximal}.
\end{defn}

\begin{lem}\label{little}
Let $\WW$ be a maximal representation of $\GG$ and let $V\in\Rep(\GG)$. If
$\Phi\in\Mor(B_V,B_\WW)$ is such that $(\Phi\tens\id)V=\WW$ then $\Phi$ is an
isomorphism.
\end{lem}

\begin{proof}
The $\cst$-algebras $B_V$ and $B_\WW$ are closures of the sets of right slices
of $V$ and $\WW$ respectively. Therefore $\Phi$ maps a dense subset of $B_V$
onto a dense subset of $B_\WW$. By the maximality of $\WW$ the map $\Phi$
increases norm:
\[
\bigl\|\Phi\bigl((\id\tens\omega)V\bigr)\bigr\|
=\bigl\|(\id\tens\omega)\WW\bigr\|
\geq\bigl\|(\id\tens\omega)V\bigr\|
\]
for any $\omega\in{A}_*$. It follows that $\Phi$ is isometric and consequently
an isomorphism.
\end{proof}

\begin{thm}\label{DefAhu}
Let $\WW$ be a maximal representation of $\GG$ and let
\begin{equation}\label{cstAhu}
\Ahu=\bigl\{(\id\tens\omega)\WW:\:\omega\in{A}_*\bigr\}^
{\scriptscriptstyle\|\cdot\|\text{\rm-closure}}.
\end{equation}
Then
\begin{enumerate}
\item\label{p1} $\Ahu$ is a non degenerate separable $\cst$-subalgebra of
$\B(H_\WW)$ and $\WW\in\M\bigl(\Ahu\tens{A}\bigr)$.
\item\label{p2} For any representation $U$ of $\GG$ there exists a unique
$\Phi_U\in\Mor\bigl(\Ahu,B_U\bigr)$ such that
\begin{equation}\label{piWU}
(\Phi_U\tens\id)\WW=U.
\end{equation}
\item\label{p3} For any pair $(B,\UU)$ such that $\UU\in\Rep(\GG)$ and $B$ is a
non degenerate $\cst$-subalgebra of $\B(H_{\UU})$ such that
$\UU\in\M(B\tens{A})$ and such that for any $U\in\Rep(\GG)$ there exists a
unique $\Phi_U\in\Mor(B,B_U)$ such that $(\Phi_U\tens\id)\UU=U$, there exists an
isomorphism $\Psi\in\Mor\bigl(\Ahu,B\bigr)$ such that $(\Psi\tens\id)\WW=\UU$.
\end{enumerate}
\end{thm}

\begin{proof}
Clearly we have $\Ahu=B_\WW$ and so Statement \eqref{p1} is just a
reformulation of the remarks at the beginning of Subsection \ref{Algs}
(cf.~\cite[Theorem 1.6]{mu}).

{\sc Ad \eqref{p2}.} To prove existence of $\Phi_U$ notice that both $U$ and
$\WW$ are subrepresentations of $U\oplus\WW$. Therefore, by Proposition
\ref{UeqV2} \eqref{UeqV20} there exist $\Phi_1\in\Mor(B_{U\oplus\WW},B_\WW)$
and $\Phi_2\in\Mor(B_{U\oplus\WW},B_U)$ such that
\[
\WW=(\Phi_1\tens\id)(U\oplus\WW),\qquad
U=(\Phi_2\tens\id)(U\oplus\WW).
\]
By Lemma \ref{little} the morphism $\Phi_1$ is an isomorphism. Then
\[
U=\bigl([\Phi_2\comp\Phi_1^{-1}]\tens\id\bigr)\WW.
\]
We let $\Phi_U=\Phi_2\comp\Phi_1^{-1}$.

Uniqueness of $\Phi_U$ is clear: if $\Phi'\in\Mor\bigl(\Ahu,B_U)$ satisfies
\[
(\Phi'\tens\id)\WW=U
\]
then for any $\omega\in{A}_*$ we have
\[
\Phi'\bigl((\id\tens\omega)\WW\bigr)=(\id\tens\omega)U=
\Phi_U\bigl((\id\tens\omega)\WW\bigr).
\]
Thus, in view of \eqref{cstAhu} we have $\Phi'=\Phi_U$.

{\sc Ad \eqref{p3}.} This is a standard consequence of the universal property
of $\bigl(\Ahu,\WW\bigr)$.
\end{proof}

\begin{rem}\label{remlit}
Note that the unique morphism $\Phi_U$ described in Theorem \ref{DefAhu}
\eqref{p2} is a surjection onto $B_U$. In particular its image does not contain
multipliers of $B_U$ which are not in $B_U$ (cf.~the proof of Lemma
\ref{little}).
\end{rem}

\begin{prop}\label{toD}
Let $\WW$ be a maximal representation of $\GG$. Then for any $\cst$-algebra $D$
and any unitary $U\in\M(D\tens{A})$ such that $(\id\tens\Delta)U=U_{12}U_{13}$
there exists a unique $\Phi_U\in\Mor\bigl(\Ahu,D\bigr)$ such that
$(\Phi_U\tens\id)\WW=U$.
\end{prop}

\begin{proof}
We can assume that the $\cst$-algebra $D$ is faithfully and non degenerately
represented on a Hilbert space $H_U$. Then $U\in\M\bigl(\K(H_U)\tens{A}\bigr)$
is a representation of $\GG$ and by Theorem \ref{DefAhu} \eqref{p2} there
exists a unique $\pi_U\in\Mor\bigl(\Ahu,B_U\bigr)$ such that \eqref{piWU}
holds.
Since $B_U$ is generated by $U$ and $U\in\M(D\tens{A})$, the identity map is a
morphism from $B_U$ to $D$. Moreover this is the only morphism from $B_U$ to
$D$ which leaves $U$ unchanged.
\end{proof}

\begin{cor}\label{RepCor}
For any $U\in\Rep(\GG)$ there exists a unique non degenerate representation
$\pi_U$ of $\Ahu$ on the Hilbert space $H_U$ such that
\[
(\pi_U\tens\id)\WW=U.
\]
Moreover the association $U\leftrightarrow\pi_U$ establishes a bijective
correspondence between $\Rep(\GG)$ and the class of all non degenerate
representations of the $\cst$-algebra $\Ahu$.
\end{cor}

\begin{proof}
Let $U\in\Rep(\GG)$. Then $U\in\M\bigl(\K(H_U)\tens{A}\bigr)$ and by
Proposition \ref{toD} there exists a unique
$\Phi_U\in\Mor\bigl(\Ahu,\K(H_U)\bigr)$ such that we have \eqref{piWU}. We let
$\pi_U$ be $\Phi_U$ considered as a map from $\Ahu$ to $\B(H_U)$. Of
course $\pi_U$ is a non degenerate representation of $\Ahu$.

Conversely, for any non degenerate representation $\pi$ of $\Ahu$ on a Hilbert
space $H$, the unitary element
$U=(\pi\tens\id)\WW\in\M\bigl(\K(H)\tens{A}\bigr)$ is a strongly continuous
unitary representation of $\GG$.
\end{proof}

\begin{rem}
Let us note that the correspondence between representations of $\GG$ and
representations of $\Ahu$ described in Corollary \ref{RepCor} is a functor
from the $\wst$-category $\Rep(\GG)$ to the $\wst$-category of
non degenerate representations of $\Ahu$. In fact it is an equivalence of
categories preserving direct sums and tensor products.
\end{rem}

\begin{defn}
The $\cst$-algebra $\Ahu$ defined in Theorem \ref{DefAhu} is called the
\emph{universal quantum group $\cst$-algebra} of $\GG$. The representation
$\WW\in\M\bigl(\Ahu\tens{A}\bigr)$ is called the \emph{universal representation}
of $\GG$.
\end{defn}

\subsection{The universal dual}

\begin{prop}\label{FirstCor}
Let $\Ahu$ be the universal quantum group $\cst$-algebra of $\GG$ and let
$\WW\in\M\bigl(\Ahu\tens{A}\bigr)$ be the universal representation. Then
\begin{enumerate}
\item\label{l1} There exists a unique
$\Delhu\in\Mor\bigl(\Ahu,\Ahu\tens\Ahu\bigr)$ such that
\begin{equation}\label{DelhWW}
\bigl(\Delhu\tens\id\bigr)\WW=\WW_{23}\WW_{13}.
\end{equation}
The morphism $\Delhu$ is coassociative and
\begin{equation}\label{linden}
\begin{array}{l@{\smallskip}}
\bigl\{\Delhu(x)(I_{\Ahu}\tens{y}):\:x,y\in\Ahu\bigr\},\\
\bigl\{(x\tens{I}_{\Ahu})\Delhu(y):\:x,y\in\Ahu\bigr\}
\end{array}
\end{equation}
are linearly dense subsets of $\Ahu\tens\Ahu$.
\item\label{l2} There exists a unique $\ehu\in\Mor\bigl(\Ahu,\CC\bigr)$ such
that
\begin{equation}\label{cou}
(\ehu\tens\id)\WW=I_A.
\end{equation}
The morphism $\ehu$ has the following property:
\begin{equation}\label{cou2}
(\id\tens\ehu)\comp\Delhu=(\ehu\tens\id)\comp\Delhu=\id.
\end{equation}
\end{enumerate}
\end{prop}

\begin{proof}
The unitary element
\[
\WW_{23}\WW_{13}\in\M\bigl(\Ahu\tens\Ahu\tens{A}\bigr)\subset
\M\bigl(\K(H_\WW)\tens\K(H_\WW)\tens{A}\bigr)=
\M\bigl(\K(H_\WW\tens{H_\WW})\tens{A}\bigr)
\]
is a strongly continuous unitary representation of $\GG$. Moreover
$\WW_{23}\WW_{13}$ is a quantum family of elements affiliated with
$B_{\WW_{23}\WW_{13}}$ generating this algebra. Therefore the identity map is a
morphism from $B_{\WW_{23}\WW_{13}}$ to $\Ahu\tens\Ahu$.

By the universal property of $\bigl(\Ahu,\WW\bigr)$, there exists a unique
$\Phi\in\Mor\bigl(\Ahu,B_{\WW_{23}\WW_{13}}\bigr)$ such that
\[
(\Phi\tens\id)\WW=\WW_{23}\WW_{13}.
\]
Let $\Delhu$ be the composition of $\Phi$ with the identity on
$B_{\WW_{23}\WW_{13}}$ considered as a morphism from $B_{\WW_{23}\WW_{13}}$ to
$\Ahu\tens\Ahu$. Then $\Delhu\in\Mor\bigl(\Ahu,\Ahu\tens\Ahu\bigr)$ satisfies
\eqref{DelhWW}.

To obtain coassociativity of $\Delhu$ we compute:
\[
\begin{array}{r@{\;=\;}l@{\smallskip}}
\Bigl(\Bigl[\bigl(\Delhu\tens\id\bigr)\comp\Delhu\Bigr]\tens\id\Bigr)\WW
&\bigl(\Delhu\tens\id\tens\id\bigr)\bigl(\Delhu\tens\id\bigr)\WW\\
&\bigl(\Delhu\tens\id\tens\id\bigr)(\WW_{23}\WW_{13})\\
&\WW_{34}\WW_{24}\WW_{14}\\
&\bigl(\id\tens\Delhu\tens\id\bigr)(\WW_{23}\WW_{13})\\
&\bigl(\id\tens\Delhu\tens\id\bigr)\bigl(\Delhu\tens\id\bigr)\WW\\
&\Bigl(\Bigl[\bigl(\id\tens\Delhu\bigr)\comp\Delhu\Bigr]\tens\id\Bigr)\WW.
\end{array}
\]
Now we can take right slice with any $\omega\in{A}_*$ and coassociativity of
$\Delhu$ follows. The the fact that the sets \eqref{linden} are contained in
$\Ahu\tens\Ahu$ and their linear density of in $\Ahu\tens\Ahu$ is proved in
the same way as \cite[Proposition 5.1]{mu} (the crucial ingredient being
\eqref{DelhWW}).

{\sc Ad \eqref{l2}.} Take $U=1\tens{I_A}\in\M\bigl(\K(\CC)\tens{A}\bigr)$. Then
$U$ is a strongly continuous unitary representation of $\GG$ and by the
universal
property of $\bigl(\Ahu,\WW\bigr)$, there exists a unique
$\ehu\in\Mor\bigl(\Ahu,\CC\bigr)$ satisfying \eqref{cou}. Notice that it follows
from \eqref{DelhWW} and \eqref{cou} that
\[
\begin{array}{r@{\;=\;}l@{\smallskip}}
\Bigl(\bigl[(\id\tens\ehu)\Delhu\bigr]\tens\id\Bigr)\WW
&(\id\tens\ehu\tens\id)\bigl(\Delhu\tens\id\bigr)\WW\\
&(\id\tens\ehu\tens\id)(\WW_{23}\WW_{13})\\
&\Bigl(I_{\Ahu}\tens\bigl[(\ehu\tens\id)\WW\bigr]\Bigr)\WW=\WW.
\end{array}
\]
In particular, for any $\omega\in{A}_*$ we obtain
\[
(\ehu\tens\id)\Delhu\bigl((\id\tens\omega)\WW\bigr)=(\id\tens\omega)
\bigl(\bigl[(\id\tens\ehu)\Delhu\bigr]\tens\id\bigr)\WW=
(\id\tens\omega)\WW
\]
and the first part of \eqref{cou2} follows. The second part is proved
analogously.
\end{proof}

\begin{defn}
Let $\Ahu$ be the universal quantum group $\cst$-algebra of $\GG$ and let
$\Delhu$ be the morphism defined in Proposition \ref{FirstCor} \eqref{l1}.
The pair $\bigl(\Ahu,\Delhu\bigr)$ will be called the \emph{universal dual} of
$\GG$.
\end{defn}

\begin{rem}
The universal dual of a quantum group is not, in general, a quantum group.
Nevertheless, as we will see, it retains a lot of structure, such as the
coinverse, scaling group and unitary coinverse.
\end{rem}

\begin{prop}\label{LamProp}
Let $W\in\M\bigl(\Hat{A}\tens{A}\bigr)$ be the reduced bicharacter for $\GG$
and $\Hat{\GG}=\bigl(\Hat{A},\Hat{\Delta}\bigr)$. Then there exists a unique
$\Hat{\Lambda}\in\Mor\bigl(\Ahu,\Hat{A}\bigr)$ such that
\[
\bigl(\Hat{\Lambda}\tens\id\bigr)\WW=W.
\]
The morphism $\Hat{\Lambda}$ satisfies
\begin{equation}\label{LamLam}
\bigl(\Hat{\Lambda}\tens\Hat{\Lambda}\bigr)\comp\Delhu
=\Hat{\Delta}\comp\Hat{\Lambda}.
\end{equation}
\end{prop}

\begin{proof}
Existence and uniqueness of $\Hat{\Lambda}$ follows from Theorem \ref{DefAhu}
\eqref{p2}.

To prove property \eqref{LamLam}, notice first  that
\[
\begin{array}{r@{\;=\;}l@{\smallskip}}
\bigl(\Hat{\Delta}\tens\id\bigr)
\bigl(\Hat{\Lambda}\tens\id\bigl)\WW
&\bigl(\Hat{\Delta}\tens\id\bigr)W=W_{23}W_{13}\\
&\bigl[\bigl(\Hat{\Lambda}\tens\id\bigr)\WW\bigr]_{23}
\bigl[\bigl(\Hat{\Lambda}\tens\id\bigr)\WW\bigr]_{13}\\
&\bigl(\Hat{\Lambda}\tens\Hat{\Lambda}\tens\id\bigr)(\WW_{23}\WW_{13})\\
&\bigl(\Hat{\Lambda}\tens\Hat{\Lambda}\tens\id\bigr)
\bigl(\Delhu\tens\id\bigr)\WW.
\end{array}
\]
Therefore for $\omega\in{A}_*$ we have
\[
\begin{array}{r@{\;=\;}l@{\smallskip}}
\Hat{\Delta}\Bigl(\Hat{\Lambda}\bigl((\id\tens\omega)\WW\bigr)\Bigr)
&(\id\tens\id\tens\omega)\bigl(\Hat{\Delta}\tens\id\bigr)
\bigl(\Hat{\Lambda}\tens\id\bigr)\WW\\
&(\id\tens\id\tens\omega)\bigl(\Hat{\Lambda}\tens\Hat{\Lambda}\tens\id\bigr)
\bigl(\Delhu\tens\id\bigr)\WW\\
&\bigl(\Hat{\Lambda}\tens\Hat{\Lambda}\bigr)
\Delhu\bigl((\id\tens\omega)\WW\bigr).
\end{array}
\]
\end{proof}

\begin{defn}
The unique $\Hat{\Lambda}\in\Mor\bigl(\Ahu,\Hat{A}\bigr)$ defined in
Proposition \ref{LamProp} is called the \emph{reducing morphism.}
\end{defn}

\begin{rem}
The reducing morphism clearly plays the role analogous to the regular
representation of a group $\cst$-algebra. We chose to name it differently in
order not to confuse it with the established notion of a regular representation
of a locally compact quantum group (\cite{kv}).
\end{rem}

\begin{prop}\label{Ahuniv}
Assume that there exists a character $\Hat{e}$ of $\Hat{A}$ such that
\begin{equation}\label{assu1}
(\id\tens\Hat{e})\Hat{\Delta}=\id
\end{equation}
and let $W\in\M\bigl(\Hat{A}\tens{A}\bigr)$ be the reduced bicharacter. Then
\begin{enumerate}
\item\label{A1} $(\Hat{e}\tens\id)W=I_A$.
\item\label{A2} $\ehu=\Hat{e}\comp\Hat{\Lambda}$.
\item\label{A3} $\Hat{\Lambda}$ is an isomorphism.
\end{enumerate}
\end{prop}

\begin{proof}
Statement \eqref{A1} is a consequence of
\[
\begin{split}
W&=(\id\tens\Hat{e}\tens\id)\bigl(\Hat{\Delta}\tens\id)W\\
&=(\id\tens\Hat{e}\tens\id)(W_{23}W_{13})
=\bigl(I_{\Hat{A}}\tens\bigl[(\id\tens\Hat{e})W\bigr]\bigr)W
\end{split}
\]
and the unitarity of $W$.

Once this is established, \eqref{A2} follows because
\[
\bigl(\bigl[\Hat{e}\comp\Hat{\Lambda}\bigr]\tens\id\bigr)\WW=
(\Hat{e}\tens\id)\bigl(\Hat{\Lambda}\tens\id\bigr)\WW=
(\Hat{e}\tens\id)W=I_A
\]
which is the defining property of $\ehu$.

We will now show that the reduced bicharacter is a maximal representation of
$\GG$. Then $\Hat{\Lambda}$ will be an isomorphism by Lemma \ref{little}. The
argument used here has already appeared in
\cite[page 177]{blanch}. We will use it in the version similar
to that of \cite[page 875]{betu}

The first observation is that it follows from the formula in Statement
\eqref{A1} that for any $\omega\in{A}_*$ we have
\begin{equation}\label{Wbig}
\bigl|\omega(I_A)\bigr|=\bigl|\Hat{e}\bigl((\id\tens\omega)W\bigr)\bigr|\leq
\bigl\|(\id\tens\omega)W\bigr\|.
\end{equation}

Let $U\in\Rep(\GG)$ and let us realize the reduced bicharacter $W$ as a
manageable multiplicative unitary on a Hilbert space $\cH$. By
\cite[Theorem 1.7]{mu} we have $W_{23}U_{12}W_{23}^*=U_{12}U_{13}$ as elements
of $\M\bigl(\K(H_U)\tens\K(\cH)\tens{A}\bigr)$. Let
$\Hat{U}=\sigma_{\K(H_U),A}(U)^*\in\M\bigl(A\tens\K(H_U)\bigr)$, where
$\sigma_{\K(H_U),A}\in\Mor\bigl(\K(H_U)\tens{A},A\tens\K(H_U)\bigr)$ is the
flip. It follows that
\[
\Hat{U}_{23}^*W_{12}=\Hat{U}_{13}W_{12}\Hat{U}_{13}^*.
\]

For each $\omega\in{A}_*$ and $\eta\in\B(H_U)_*$ we define
$\omega_\eta\in{A}_*$ by
\[
\omega_\eta(x)=(\omega\tens\eta)\bigl(\Hat{U}^*(x\tens{I_{H_U}})\bigr).
\]
We have
\begin{equation}\label{splitting}
\begin{split}
(\id\tens\omega_\eta)W&=(\id\tens\omega\tens\eta)
\bigl(\Hat{U}_{23}^*W_{12}\bigr)\\
&=(\id\tens\omega\tens\eta)\bigl(\Hat{U}_{13}W_{12}\Hat{U}_{13}^*\bigr)\\
&=(\id\tens\eta)
\Bigl[\Hat{U}\bigl[\bigl((\id\tens\omega)W\bigr)\tens{I}\bigr]\Hat{U}^*\Bigr].
\end{split}
\end{equation}
Also
\begin{equation}\label{etaI}
\begin{split}
\eta\bigl((\id\tens\omega)U\bigr)
&=\eta\bigl((\omega\tens\id)(\Hat{U}^*)\bigr)\\
&=(\omega\tens\eta)(\Hat{U}^*)=\omega_\eta(I_A).
\end{split}
\end{equation}
Now using \eqref{etaI}, \eqref{Wbig} and \eqref{splitting} we have
\[
\begin{split}
\bigl|\eta\bigl((\id\tens\omega)U\bigr)\bigr|
&=\bigl|\omega_\eta(I_A)\bigr|\leq\bigl\|(\id\tens\omega_\eta)W\bigr\|\\
&=\Bigl\|(\id\tens\eta)
\Bigl[\Hat{U}\bigl[\bigl((\id\tens\omega)W\bigr)\tens{I}\bigr]\Hat{U}^*\Bigr]
\Bigr\|\leq\|\eta\|\bigl\|(\id\tens\omega)W\bigr\|.
\end{split}
\]
Since for any $t\in\B(H_U)$ we have
$\|t\|=
\sup\bigl\{\bigl|\eta(t)\bigr|:\:\eta\in\B(H_U)_*,\:\|\eta\|=1\bigr\}$,
we conclude that for any $\omega\in{A}_*$
\[
\bigl\|(\id\tens\omega)U\bigr\|\leq\bigl\|(\id\tens\omega)W\bigr\|.
\]
\end{proof}

\begin{rem}
The assumption \eqref{assu1} in Proposition \ref{Ahuniv} is in fact equivalent
to the formula in Statement \eqref{A1} of that Proposition (cf.~the proof of
Proposition \ref{FirstCor} \eqref{l2}).
\end{rem}

\begin{prop}\label{tauhu}
Let $\bigl(\Ahu,\WW\bigr)$ be the universal quantum group $\cst$-algebra of and
the universal representation of $\GG$. Then
\begin{enumerate}
\item\label{T1} for any $t\in\RR$ there exists a unique
$\tauhu_t\in\Mor(\Ahu,\Ahu)$ such that
\begin{equation}\label{tautau}
(\tauhu_t\tens\id)\WW=(\id\tens\tau_{-t})\WW.
\end{equation}
Moreover $(\tauhu_t)_{t\in\RR}$ is a one parameter group of automorphisms of
$\Ahu$.
\item\label{T2} For any $x\in\Ahu$ the map $\RR\ni{t}\mapsto\tauhu_t(x)\in\Ahu$
is continuous.
\item\label{T3} For any $t\in\RR$ we have
\begin{equation}\label{tauhQG}
(\tauhu_t\tens\tauhu_t)\comp\Delhu=\Delhu\comp\tauhu_t
\end{equation}
and $\ehu\comp\tauhu_t=\ehu$.
\item\label{T4} If $(\Hat{\tau}_t)_{t\in\RR}$ is the scaling group of
$\Hat{\GG}$ and $\Hat{\Lambda}\in\Mor\bigl(\Ahu,\Hat{A}\bigr)$ is the reducing
morphism then for any $t\in\RR$ we have
\begin{equation}\label{tauLam}
\Hat{\tau}_t\comp\Hat{\Lambda}=\Hat{\Lambda}\comp\tauhu_t.
\end{equation}
\end{enumerate}
\end{prop}

\begin{proof}
{\sc Ad \eqref{T1}.}
For any $t\in\RR$ the element $(\id\tens\tau_{-t})\WW\in
\M\bigl(\K(H_\WW)\tens{A}\bigr)$ is a representation of $\GG$. Therefore, by
the universal property of $\bigl(\Ahu,\WW\bigr)$, there exists a unique
$\tauhu_t\in\Mor\bigl(\Ahu,\Ahu\bigr)$ such that \eqref{tautau} holds. It is
easy to see that $(\tauhu_t)_{t\in\RR}$ is a one parameter group of
automorphisms of $\Ahu$.

{\sc Ad \eqref{T2}.} Take $x\in\Ahu$ and $\eps>0$. There exists a functional
$\omega\in{A}_*$ such that
\[
\|x-(\id\tens\omega)\WW\|<\textstyle{\frac{\eps}{3}}
\]
and the map $\RR\ni{t}\mapsto\omega\comp\tau_{-t}$ is norm continuous.
Therefore for $t,s\in\RR$
\[
\begin{array}{r@{\;}c@{\;}l@{\medskip}}
\bigl\|\tauhu_t(x)-\tauhu_s(x)\bigr\|&\leq&
\bigl\|\tauhu_t(x)-\tauhu_t\bigl((\id\tens\omega)\WW\bigr)\bigr\|\\
&&+\bigl\|\tauhu_t\bigl((\id\tens\omega)\WW\bigr)
-\tauhu_s\bigl((\id\tens\omega)\WW\bigr)\bigr\|\\
&&+\bigl\|\tauhu_s\bigl((\id\tens\omega)\WW\bigr)-\tauhu_s(x)\bigr\|.
\end{array}
\]
The first and third terms are each smaller than $\frac{\eps}{3}$
and the middle term
\[
\begin{array}{r@{\;}c@{\;}l@{\medskip}}
\bigl\|\tauhu_t\bigl((\id\tens\omega)\WW\bigr)
-\tauhu_s\bigl((\id\tens\omega)\WW\bigr)\bigr\|
&=&
\bigl\|(\id\tens\omega)\bigl((\tauhu_t\tens\id)\WW-
(\tauhu_s\tens\id)\WW\bigr)\bigr\|\\
&=&
\bigl\|(\id\tens\omega)\bigl((\id\tens\tau_{-t})\WW-
(\id\tens\tau_{-s})\WW\bigr)\bigr\|\\
&=&\bigl\|\bigl(\id\tens
[\omega\comp\tau_{-t}-\omega\comp\tau_{-s}]\bigr)\WW\bigr\|\\
&\leq&\bigl\|\omega\comp\tau_{-t}-\omega\comp\tau_{-s}\bigr\|
\end{array}
\]
is smaller than $\frac{\eps}{3}$ for $s$ sufficiently close to $t$.

{\sc Ad \eqref{T3}.} Just as in the proof of formula \eqref{LamLam}
(Proposition \ref{LamProp}) we can easily show that
\[
\Bigl(\bigl(\Delhu\comp\tauhu_t\bigr)\tens\id\Bigr)\WW=
\Bigl(\bigl[(\tauhu_t\tens\tauhu_t)\comp\Delhu\bigr]\tens\id\Bigr)\WW
\]
and \eqref{tauhQG} follows. Similarly for any $t\in\RR$ and $\omega\in{A}_*$
\[
\begin{array}{r@{\;=\;}l@{\smallskip}}
\ehu\Bigl(\tauhu_t\bigl((\id\tens\omega)\WW\bigr)\Bigr)
&\ehu\bigl((\id\tens\omega)(\tauhu_t\tens\id)\WW\bigr)\\
&\ehu\bigl((\id\tens\omega)(\id\tens\tau_{-t})\WW\bigr)\\
&\omega\Bigl(\tau_{-t}\bigl((\ehu\tens\id)\WW\bigr)\Bigr)\\
&\omega\bigl(\tau_{-t}(I_A)\bigr)=\omega(I_A)\\
&\omega\bigl((\ehu\tens\id)\WW\bigr)\\
&\ehu\bigl((\id\tens\omega)\WW\bigr)
\end{array}
\]
proves the other formula.

{\sc Ad \eqref{T4}.} Let $W$ be the reduced bicharacter for
$\bigl(\GG,\Hat{\GG})$.
First let us see that for any $t\in\RR$ we have
\begin{equation}\label{tmt}
(\Hat{\tau}_t\tens\id)W=(\id\tens\tau_{-t})W.
\end{equation}
Indeed, we can realize $W$ as a modular multiplicative unitary on some Hilbert
space $\cH$. Then (\cite{mu,mmu}) the scaling group of $\GG$  is given by
$\tau_t(a)=Q^{2it}aQ^{-2it}$ where $Q$ is one of the two positive self
adjoint operators appearing in the definition of modularity of $W$. Similarly
the scaling group of $\Hat{\GG}$ is
$\Hat{\tau}_t(x)=\Hat{Q}^{2it}x\Hat{Q}^{-2it}$, where $\Hat{Q}$ is the other
positive self adjoint operator. Now since $W$ commutes with $\Hat{Q}\tens{Q}$,
we obtain \eqref{tmt}.

Now for $t\in\RR$ we can compute
\[
\begin{array}{r@{\;=\;}l@{\smallskip}}
\bigl(\Hat{\Lambda}\tens\id\bigr)(\tauhu_t\tens\id)\WW
&\bigl(\Hat{\Lambda}\tens\id\bigr)(\id\tens\tau_{-t})\WW\\
&(\id\tens\tau_{-t})\bigl(\Hat{\Lambda}\tens\id\bigr)\WW\\
&(\id\tens\tau_{-t})W=(\Hat{\tau}_t\tens\id)W\\
&(\tauhu_t\tens\id)\bigl(\Hat{\Lambda}\tens\id\bigr)\WW.
\end{array}
\]
As before, the resulting formula
\[
\bigl(\bigl[\Hat{\Lambda}\comp\tauhu_t\bigr]\tens\id\bigr)\WW
=\bigl(\bigl[\Hat{\tau}_t\comp\Hat{\Lambda}\bigr]\tens\id\bigr)\WW
\]
suffices to have \eqref{tauLam}.
\end{proof}

\begin{lem}\label{WRR}
Let $W\in\M\bigl(\Hat{A}\tens{A}\bigr)$ be the reduced bicharacter. Then
\[
W^{\Hat{R}\tens{R}}=W.
\]
\end{lem}

\begin{proof}
Let $\omega\in\Hat{A}_*$ and $\mu\in{A}_*$ be analytic for
(the transpose of) $(\Hat{\tau}_t)$ and $(\tau_t)$ respectively. Using the fact
that for any $\nu\in\Hat{A}_*$ and $\lambda\in{A}_*$
\[
\begin{split}
\kappa\bigl((\nu\tens\id)W\bigr)&=(\nu\tens\id)(W^*),\\
\widehat{\kappa}\bigl((\lambda\tens\id)\widehat{W}\bigr)&
=(\lambda\tens\id)\bigl(\widehat{W}^*\bigr)\\
\end{split}
\]
we have
\[
\begin{split}
\bigl([\omega\comp\Hat{\tau}_{\frac{i}{2}}]\tens
[\mu\comp\tau_{\frac{i}{2}}]\bigr)\bigl(W^{\Hat{R}\tens{R}}\bigr)
&=\bigl(\bigl[\omega\comp\Hat{\tau}_{\frac{i}{2}}\comp\Hat{R}]\tens
[\mu\comp\tau_{\frac{i}{2}}\comp{R}]\bigr)W\\
&=\bigl(\bigl[\omega\comp\Hat{\tau}_{\frac{i}{2}}
\comp\Hat{R}]\tens\mu\bigr)(W^*)\\
&=\bigl(\mu\tens\bigl[\omega\comp\Hat{\tau}_{\frac{i}{2}}
\comp\Hat{R}]\bigr)\Hat{W}\\
&=(\mu\tens\omega)\bigl(\Hat{W}^*\bigr)=(\omega\tens\mu)W\\
&=\bigl([\omega\comp\Hat{\tau}_{\frac{i}{2}}]\tens
[\mu\comp\tau_{\frac{i}{2}}]\bigr)W,
\end{split}
\]
where in the last step we used holomorphic continuation of
\[
\bigl(\Hat{\tau}_t\tens\tau_t\bigr)W=W
\]
which follows from \eqref{tmt}. The conclusion follows from the fact that
functionals of the form $\bigl([\omega\comp\Hat{\tau}_{\frac{i}{2}}]\tens
[\mu\comp\tau_{\frac{i}{2}}]\bigr)$ separate points of $\Hat{A}\tens{A}$.
\end{proof}

Before the next proposition let us state a remark which will come handy in the
proof.

\begin{rem}
Let $C$ be a $\cst$-algebra and let $S$ be an anti-morphism from $\Ahu$ to
$C$, i.e.~$S\in\Mor\bigl(\Ahu,C^{\mathrm{op}}\bigr)$. Then for any
$\omega\in{A}_*$ we have
\begin{equation}\label{SRw}
(\id\tens\omega)\comp(S\tens{R})=S\comp\bigl(\id\tens[\omega\comp{R}]\bigr).
\end{equation}
\end{rem}

\begin{prop}
Let $\bigl(\Ahu,\WW\bigr)$ be the universal quantum group $\cst$-algebra and the
universal representation of $\GG$. Then
\begin{enumerate}
\item\label{R1} There exists a unique anti-automorphism $\Rhu$ of $\Ahu$ such
that for any $\omega\in{A}_*$
\begin{equation}\label{DefRhu}
\Rhu\bigl((\id\tens\omega)\WW\bigr)=\bigl((\id\tens\omega)\WW^\cc\bigr)^\top.
\end{equation}
\item\label{R2} $\Rhu$ is involutive.
\item\label{R3} Let $\Hat{\sigma}$ be the flip on $\Ahu\tens\Ahu$. Then
\begin{equation}\label{first}
\Hat{\sigma}\comp\bigl(\Rhu\tens\Rhu\bigr)\comp\Delhu=\Delhu\comp\Rhu.
\end{equation}
Moreover if $(\tauhu_t)_{t\in\RR}$ is the one parameter group of automorphisms
of $\Ahu$ defined in Proposition \ref{tauhu} then for any $t\in\RR$ we have
\begin{equation}\label{second}
\tauhu_t\comp\Rhu=\Rhu\comp\tauhu_t.
\end{equation}
\item\label{R4} If $\Hat{R}$ is the unitary coinverse of $\Hat{\GG}$ and
$\Hat{\Lambda}\in\Mor\bigl(\Ahu,\Hat{A}\bigr)$ is the reducing morphism then
we have
\[
\Hat{R}\comp\Hat{\Lambda}=\Hat{\Lambda}\comp\Rhu.
\]
\end{enumerate}
\end{prop}

\begin{proof}
{\sc Ad \eqref{R1}.} The contragradient representation $\WW^\cc$ of $\WW$ is a
strongly continuous unitary representation of $\GG$. Therefore, by the
universal property of $\bigl(\Ahu,\WW\bigr)$, there exists a
unique morphism $\theta\in\Mor\bigl(\Ahu,B_{\WW^\cc}\bigr)$ such that
\begin{equation}\label{thetWc}
(\theta\tens\id)\WW=\WW^\cc
\end{equation}
(cf.~Theorem \ref{DefAhu} \eqref{p2}). Clearly
$B_{\WW^\cc}=\bigl(\Ahu\bigr)^\top$ and we can define a map
\[
\Rhu:\Ahu\ni{x}\longmapsto\theta(x)^\top\in\Ahu.
\]
It is easy to see that so defined $\Rhu$ is an anti-morphism of $\Ahu$ to
itself which satisfies \eqref{DefRhu} which determines this anti-morphism
uniquely. The map $\Rhu$ is an anti-automorphism of $\Ahu$. This follows for
example from the fact that $\bigl(\Rhu\bigr)^2=\id$ established below.

{\sc Ad \eqref{R2}.} Let us take contragradient representations of both sides
of \eqref{thetWc}. Then
\begin{equation}\label{RRW}
\bigl(\Rhu\tens{R}\bigr)\WW=\WW.
\end{equation}
Now applying $\bigl(\Rhu\tens{R}\bigr)$ to both sides of \eqref{RRW} and
then using this equation we arrive at
\[
\bigl(\bigl(\Rhu\bigr)^2\tens\id\bigr)\WW=\bigl(\Rhu\tens{R}\bigr)\WW=\WW.
\]
Thus for any $\omega\in{A}_*$ we have
\[
\bigl(\Rhu\bigr)^2\bigl((\id\tens\omega)\WW\bigr)=(\id\tens\omega)
\bigl(\bigl(\Rhu\bigr)^2\tens\id\bigr)\WW=(\id\tens\omega)\WW
\]
and it follows $\bigl(\Rhu\bigr)^2=\id$.

{\sc Ad \eqref{R3}.} Let us begin with \eqref{second}. For any $t\in\RR$ we
have:
\[
\begin{split}
\bigl(\bigl[\Rhu\comp\tauhu_t\bigr]\tens{R}\bigr)\WW
&=\bigl(\Rhu\tens{R}\bigr)(\tauhu_t\tens\id)\WW
=\bigl(\Rhu\tens{R}\bigr)(\id\tens\tau_{-t})\WW\\
&=\bigl(\Rhu\tens[R\comp\tau_{-t}]\bigr)\WW
=\bigl(\Rhu\tens[\tau_{-t}\comp{R}]\bigr)\WW\\
&=(\id\tens\tau_{-t})\bigl(\Rhu\tens{R}\bigr)\WW
=(\id\tens\tau_{-t})\WW\\
&=(\tauhu_t\tens\id)\WW=(\tauhu_t\tens\id)\bigl(\Rhu\tens{R}\bigr)\WW\\
&=\bigl(\bigl[\tauhu_t\comp\Rhu\bigr]\tens{R}\bigr)\WW.
\end{split}
\]
Now let us take $\omega\in{A}_*$. Then with information from the above
computation and using formula \eqref{SRw} twice (once with $S=\Rhu\comp\tauhu_t$
and then with $S=\tauhu_t\comp\Rhu$) we get
\[
\begin{split}
\bigl(\Rhu\comp\tauhu_t\bigr)
\Bigl(\bigl(\id\tens[\omega\comp{R}]\bigr)\WW\Bigr)
&=\Rhu\Bigl(\bigl(\id\tens[\omega\comp{R}]\bigr)
\bigl((\tauhu_t\tens\id)\WW\bigr)\Bigr)\\
&=(\id\tens\omega)\bigl(\bigl[\Rhu\comp\tauhu_t\bigr]\tens{R}\bigr)\WW\\
&=(\id\tens\omega)\bigl(\bigl[\tauhu_t\comp\Rhu\bigr]\tens{R}\bigr)\WW\\
&=\bigl(\tauhu_t\comp\Rhu\bigr)
\Bigl(\bigl(\id\tens[\omega\comp{R}]\bigr)\WW\Bigr)
\end{split}
\]
Now since $\bigl\{\omega\comp{R}:\:\omega\in{A}_*\bigr\}=A_*$ and $\Ahu$ is
the closure of the set of right slices of $\WW$ we get \eqref{second}.

For the proof of \eqref{first} let us first notice that the defining
property of $\Delhu$ and \eqref{RRW} imply that
\[
\begin{split}
\bigl(\Rhu\tens\Rhu\tens{R}\bigr)\bigl(\Delhu\tens\id\bigr)\WW
&=\bigl(\Rhu\tens\Rhu\tens{R}\bigr)(\WW_{23}\WW_{13})\\
&=\bigl[\bigl(\Rhu\tens{R}\bigr)\WW\bigr]_{13}
\bigl[\bigl(\Rhu\tens{R}\bigr)\WW\bigr]_{23}\\
&=\WW_{13}\WW_{23}.
\end{split}
\]
Therefore
\[
\begin{split}
\Bigl(\bigl[\Hat{\sigma}\comp\bigl(\Rhu\tens\Rhu\bigr)\comp\Delhu\bigr]
\tens{R}\Bigr)\WW&=\WW_{23}\WW_{13}
=\bigl(\Delhu\tens\id\bigr)\WW\\
&=\bigl(\Delhu\tens\id\bigr)\bigl(\Rhu\tens\Rhu\bigr)\WW\\
&=\bigl(\bigl[\Delhu\comp\Rhu]\tens{R}\bigr)\WW
\end{split}
\]
and in the same way as in the proof of \eqref{second} we obtain \eqref{first}.
This time formula \eqref{SRw} must also be used twice: once with
$S=\Hat{\sigma}\comp\bigl(\Rhu\tens\Rhu\bigr)\comp\Delhu$ and then with
$S=\Delhu\comp\Rhu$.

{\sc Ad \eqref{R4}.} Let $W\in\M\bigl(\Hat{A}\tens{A}\bigr)$ be the reduced
bicharacter. Using \eqref{RRW} in the second and Lemma \ref{WRR} in the fourth
step we have
\[
\begin{split}
\bigl(\bigl[\Hat{\Lambda}\comp\Rhu\bigr]\tens{R}\bigr)\WW
&=\bigl(\Hat{\Lambda}\tens\id\bigr)\bigl(\Rhu\tens{R}\bigr)\WW\\
&=\bigl(\Hat{\Lambda}\tens\id\bigr)\WW=W=W^{\Hat{R}\tens{R}}\\
&=\bigl(\Rhu\tens{R}\bigr)\bigl(\Hat{\Lambda}\tens\id\bigr)\WW
=\bigl(\bigl[\Hat{R}\comp\Hat{\Lambda}\bigr]\tens{R}\bigr)\WW.
\end{split}
\]
Again, as in proofs of \eqref{R2} and \eqref{R3} we can use \eqref{SRw} once
with $S=\Hat{\Lambda}\comp\Rhu$ and then with $S=\Hat{R}\comp\Hat{\Lambda}$
and appeal to the fact that $R$ is a homeomorphism for the ultraweak topology
on $A$.
\end{proof}

\begin{prop}
Let $\bigl(\Ahu,\WW\bigr)$ be the universal quantum group $\cst$-algebra and the
universal representation of $\GG$. Then there exists a unique closed linear
operator $\kaphu$ on the Banach space $\Ahu$ such that
\begin{equation}\label{core}
\bigl\{(\id\tens\omega)(\WW^*):\:\omega\in{A}_*\bigr\}
\end{equation}
is a core for $\kaphu$ and
\[
\kaphu\bigl((\id\tens\omega)(\WW^*)\bigr)=(\id\tens\omega)\WW.
\]
Moreover
\begin{enumerate}
\item\label{kap1} the domain of $\kaphu$ is an algebra and $\kaphu$ is
anti-multiplicative: $\kaphu(xy)=\kaphu(y)\kaphu(x)$ for all
$x,y\in\Dom(\kaphu)$.
\item\label{kap2} For any $x\in\Dom(\kaphu)$ the element $\kaphu(x)^*$ belongs
to $\Dom(\kaphu)$ and we have $\kaphu\bigl(\kaphu(x)^*\bigr)^*=x$.
\item\label{kap3} $\kaphu=\Rhu\comp\tauhu_{\frac{i}{2}}$.
\end{enumerate}
\end{prop}

\begin{proof}
Take $\eta\in\B(H_\WW)_*$ and $\omega\in{A}_*$. Then with repeated use of a
variation of formula \eqref{SRw} we compute:
\[
\begin{split}
\bigl(\eta\comp\Rhu\comp\tauhu_t\bigr)
\Bigl[\bigl(\id\tens[\omega\comp{R}]\bigr)(\WW^*)\Bigr]
&=\bigl(\bigl[\eta\comp\Rhu\comp\tauhu_t\bigr]
\tens[\omega\comp{R}]\bigr)(\WW^*)\\
&=\bigl(\bigl[\eta\comp\Rhu\bigr]\tens[\omega\comp{R}]\bigr)
\bigl[(\tauhu_t\tens\id)(\WW^*)\bigr]\\
&=\bigl(\bigl[\eta\comp\Rhu\bigr]\tens[\omega\comp{R}]\bigr)
\bigl[(\id\tens\tau_{-t})(\WW^*)\bigr]\\
&=\bigl(\bigl[\eta\comp\Rhu\bigr]
\tens[\omega\comp{R}\comp\tau_{-t}]\bigr)(\WW^*)\\
&=(\omega\comp{R}\comp\tau_{-t})
\Bigl[\bigl(\bigl[\eta\comp\Rhu\bigr]\tens\id\bigr)(\WW^*)\Bigr]\\
&=\omega\Bigl[(R\comp\tau_{-t})
\bigl(\bigl[\eta\comp\Rhu\bigr]\tens\id\bigr)(\WW^*)\Bigr].
\end{split}
\]
We shall now take the limit $t\to\frac{i}{2}$. By Proposition \ref{kapU} and
properties of analytic generators of one parameter groups, the last term above
converges to
\[
\begin{split}
\omega\Bigl[\bigl(\bigl[\eta\comp\Rhu\bigr]\tens\id\bigr)\WW\Bigr]
&=(\omega\comp{R})
\Bigl[\bigl(\bigl[\eta\comp\Rhu\bigr]\tens{R}\bigr)\WW\Bigr]\\
&=\bigl(\eta\tens[\omega\comp{R}]\bigr)
\Bigl[\bigl(\Rhu\tens{R}\bigr)\WW\Bigr]
=\bigl(\eta\tens[\omega\comp{R}]\bigr)\WW
\end{split}
\]
(cf.~\eqref{RRW}). Since functionals of the form $\omega\comp{R}$ fill up all
of $A_*$, we find that for any $\omega\in{A}_*$ the element
$(\id\tens\omega)(\WW^*)$ belongs to the domain of the analytic extension of
the group $(\tauhu_t)_{t\in\RR}$ to the point $t=\frac{i}{2}$ and we have
\[
\bigl(\Rhu\comp\tauhu_{\frac{i}{2}}\bigr)\bigl((\id\tens\omega)(\WW^*)\bigr)
=(\id\tens\omega)\WW.
\]
Moreover the set \eqref{core} is dense in $\Ahu$ and is
$(\tauhu_t)_{t\in\RR}$-invariant. Indeed, recall that all automorphisms
$(\tau_t)_{t\in\RR}$ are ultraweakly continuous, so
\[
\tauhu_t\bigl((\id\tens\omega)(\WW^*)\bigr)=
(\id\tens\omega)(\tauhu_t\tens\id)(\WW^*)=
(\id\tens\omega)(\id\tens\tau_{-t})(\WW^*)
=\bigl(\id\tens[\omega\comp\tau_{-t}]\bigr)(\WW^*)
\]
belongs to \eqref{core}. It is a simple observation (cf.~e.g.~\cite[Proposition
F.5]{mnw}) that this implies that \eqref{core} must be a core for
$\tauhu_{\frac{i}{2}}$.
If we now put $\kaphu=\Rhu\comp\tauhu_{\frac{i}{2}}$ then $\kaphu$ is a closed
operator on $\Ahu$. Clearly there at most one operator with a given core and
prescribed action on this core. Of course, w have \eqref{kap3}. Properties
\eqref{kap1}
and \eqref{kap2} are well known facts from the theory analytic extensions of
one parameter groups of automorphisms (\cite{zsido}).
\end{proof}


\begin{thebibliography}{66}
\bibitem{bs}
{\sc S.~Baaj \& G.~Skandalis,} Unitaries muliplicatifs et dualit\'{e}
pour les poiduits crois\'{e}s de $\mathrm{C}^*$-alg\'{e}bres.
\emph{Ann.~Scient.~\'{E}c.~Norm.~Sup.} $4^{\rm e}$ s\'{e}rie, t.~\textbf{26}
(1993), 425--488.
\bibitem{bsv}
{\sc S.~Baaj, G.~Skandalis \& S.~Vaes,} Non-semi-regular quantum groups coming
from number theory. \emph{Commun.~Math.~Phys.} \textbf{235} no.~1 (2003),
139--167.
\bibitem{betu}
{\sc E.~Bedos \& L.~Tuset:} Amenability and co-amenability for locally compact
quantum groups. \emph{Int.~J.~Math.} \textbf{14} no.~8 (2003), 865--884.
\bibitem{blanch}
{\sc E.~Blanchard:} D\'eformations de $\mathrm{C}^*$-alg\`ebres de Hopf.
\emph{Bull.~Soc.~Math.~France} \textbf{124} (1996), 141--215.
\bibitem{dix}
{\sc J.~Dixmier:} \emph{$\mathrm{C}^*$-algebras and their representations,}
North Holland 1977.
\bibitem{Wst}
{\sc P.~Ghez, R.~Lima \& J.E.~Roberts:} $\mathrm{W}^*$-categories.
\emph{Pac.~J.~Math.} \textbf{120} no.~1 (1985), 79--109.
\bibitem{univLCQG}
{\sc J.~Kustermans:} Locally compact quantum groups in the universal setting.
\emph{Int.~J.~Math.} \textbf{12} (2001), 289-338.
\bibitem{kv}
{\sc J.~Kustermans \& S.~Vaes:} Locally compact quantum groups.
\emph{Ann.~Sci.~Ec.~Norm.~Sup.} \textbf{33} no.~4 (2000), 837--934.
\bibitem{lance}
{\sc E.~Lance:} \emph{Hilbert $\mathrm{C}^*$-modules, an operator algebraist's
toolkit.} Cambridge 1995.
\bibitem{mnw}
{\sc T.~Masuda, Y.~Nakagami \& S.L.~Woronowicz:} A $\mathrm{C}^*$-algebraic
framework for quantum groups. \emph{Int.~J.~Math.} \textbf{14} no.~9 (2003)
903--1001.
\bibitem{pu-slw}
{\sc W.~Pusz \& S.L.~Woronowicz:} On quantum group of unitary operators.
Quantum `$az+b$' group. \emph{Twenty years of Bia{\l}owie\.za}, World
Scientific Monograph Series in Mathematics, Vol.~8, World Scientific (2005),
pp.~229--261.
\bibitem{nazb}
{\sc P.M.~So{\l}tan:} New Quantum ``$az+b$'' groups. \emph{Rev.~Math.~Phys.}
\textbf{17} no.~{3} (2005), 313--364.
\bibitem{mmu}
{\sc P.M.~So{\l}tan \& S.L.~Woronowicz:} A remark on manageable
multiplicative unitaries. \emph{Lett.~Math.~Phys.} \textbf{57} (2001),
239--252.
% \bibitem{indu}
% {\sc S.~Vaes:} A new approach to induction and imprimitivity results.
% \emph{J.~Func.~Anal.} \textbf{229} (2005), 317--374.
\bibitem{VDhaar}
{\sc A.~Van Daele,} The Haar measure on some locally compact
quantum groups. \emph{Preprint K.U.~Leuven} (2001), {\tt math.OA/0109004.}
\bibitem{pseu}
{\sc S.L.~Woronowicz:} Pseudogroups, pseudospaces and Pontryagin duality.
\emph{Proceedings of the International Conference on Mathematical Physics,
Lausanne 1979} Lecture Notes in Physics, \textbf{116}.
\bibitem{pseudogr}
{\sc S.L.~Woronowicz:} Compact matrix pseudogroups. \emph{Comm.~Math.~Phys.}
\textbf{111} no.~4 (1987), 613--665.
\bibitem{TKSUN}
{\sc S.L.~Woronowicz:} Tannaka-Krein duality for compact matrix pseudogroups.
Twisted $SU(N)$ groups. \emph{Invent.~Math.} \textbf{93} (1988), 35--76.
\bibitem{unbo}
{\sc S.L.~Woronowicz:} Unbounded elements affiliated with
$\mathrm{C}^*$-algebras and non-compact quantum groups.
\emph{Commun.~Math.~Phys.} \textbf{136} (1991), 399--432.
\bibitem{gen}
{\sc S.L.~Woronowicz:} $\mathrm{C}^*$-algebras generated by unbounded elements.
\emph{Rev.~Math. Phys.} \textbf{7} no.~3 (1995), 481--521.
\bibitem{mu}
{\sc S.L.~Woronowicz:} From multiplicative unitaries to quantum groups.
\emph{Int.~J.~Math.} \textbf{7} no.~1 (1996), 127--149.
\bibitem{azb}
{\sc S.L.~Woronowicz:} Quantum `$az+b$' group on complex plane,
\emph{Int.~J.~Math.} \textbf{12}, no.~4 (2001), 461--503.
\bibitem{haar}
{\sc S.L.~Woronowicz:} Haar weight on some quantum groups. \emph{Group 24:
Physical and mathematical aspects of symmetries, Proceedings of the
24th International Colloquium on Group Theoretical Methods in Physics Paris, 15
- 20 July 2002,} Institute of Physics, Conference Series Number \textbf{173},
pp.~763--772.
\bibitem{axb}
{\sc S.L.~Woronowicz, \& S.~Zakrzewski:} Quantum `$ax+b$' group,
\emph{Rev.~Math.~Phys.} \textbf{14}, nos.~7 \& 8 (2002), 797--828.
\bibitem{zsido}
{\sc L.~Zsid\'{o}:} On the characterization of the analytic generator of
$*$-automorphism groups. In \emph{Operator algebras and applications,
Proc.~Symp.~Pure.~Math.} \textbf{38}, AMS Providence, R.I.~(1982),
pp.~381--383.
\end{thebibliography}
\end{document}